\documentclass[12pt,a4paper]{article}
\usepackage{amsmath}
\usepackage{amsfonts}
\usepackage{amsthm}
\usepackage{amssymb}

\usepackage{enumerate}
\theoremstyle{plain}
\newtheorem{theo}{Theorem}[section]
\newtheorem{prop}[theo]{Proposition}
\theoremstyle{definition}
\newtheorem{definition}[theo]{Definition}
\theoremstyle{remark}
\numberwithin{equation}{section}

\newcommand{\R}{\mathbb{R}}
\newcommand{\N}{\mathbb{N}}

\title{Doubling inequality at the boundary for the Kirchhoff-Love plate's equation with supported conditions
\thanks{ER, ES and SV are supported by 
PRIN 201758MTR2 ``Direct and inverse problems for partial differential equations: theoretical
aspects and applications''. ER and SV are also funded by Progetto GNAMPA 2019 ``Propriet\`a delle soluzioni di equazioni alle
derivate parziali e applicazioni ai problemi inversi" Istituto Nazionale di Alta Matematica (INdAM).}}
\author{Antonino Morassi\thanks{Dipartimento Politecnico di Ingegneria e Architettura,
Universit\`a degli Studi di Udine, via Cotonificio 114, 33100
Udine, Italy. E-mail: \textsf{antonino.morassi@uniud.it}} \ Edi
Rosset\thanks{Dipartimento di Matematica e Geoscienze,
Universit\`a degli Studi di Trieste, via Valerio 12/1, 34127
Trieste, Italy. E-mail: \textsf{rossedi@units.it}}  \ Eva
Sincich\thanks{Dipartimento di Matematica e Geoscienze,
Universit\`a degli Studi di Trieste, via Valerio 12/1, 34127
Trieste, Italy. E-mail: \textsf{esincich@units.it}} 
\ and
Sergio Vessella\thanks{Dipartimento di Matematica e Informatica ``Ulisse Dini'', Universit\`a degli Studi di Firenze,Viale Morgagni 67/a,
50134 Firenze, Italy. E-mail:
\textsf{sergio.vessella@unifi.it}}}

\begin{document}

\maketitle

\begin{center}
 \textbf{Abstract} 
\end{center}

\noindent In this article we derive a doubling inequality at the boundary for solutions to the Kirchhoff-Love isotropic plate's equation satisfying supported boundary conditions. To this end, we combine
the use of a suitable conformal mapping which flattens the boundary and a reflection argument which guarantees the needed regularity of the extended solution. We finally apply inequalities of Carleman type in order to derive the result.  The latter implies Strong Unique Continuation Property at the boundary (SUCPB).
\medskip


\medskip

\medskip


\section{Introduction} \label{sec:
introduction}

In this paper we are mainly concerned with a Strong Unique Continuation Property at the Boundary for the Kirchhoff-Love isotropic plate's equation satisfying the so-called supported conditions. 
\noindent We denote by $\Omega\subset \mathbb{R}^2$ the middle surface of a thin plate having uniform thickness $h$. Working in the framework of linear elasticity for infinitesimal deformations and under the assumptions of the Kirchhoff-Love theory, the transverse displacement $u$ of the plate satisfies the following fourth order partial differential equation 
\begin{equation}
    \label{eq:equazione_piastra_intro}
   {\rm div}\left ({\rm div} \left ( B(1-\nu)\nabla^2 u + B\nu \Delta u I_2 \right ) \right )=0, \qquad\hbox{in
    } \Omega ,
\end{equation}
where $B$ is the bending stiffness and $\nu$ is the Poisson's coefficient. 
\noindent We are interested in analyzing the case of supported conditions, namely 
\begin{center}
\( {\displaystyle \left\{
\begin{array}{lr}
     u=0, \quad &\hbox{on } \Gamma,\\
     B(1-\nu)\partial^2_{ij}u\,n_in_j +B\nu\Delta u=0 \quad &\hbox{on } \Gamma ,
\end{array}
\right. } \) \vskip -4.2em
\begin{eqnarray}
& & \label{eq:bdry_cond1_intro}\\
& & \label{eq:bdry_cond2_intro}
\end{eqnarray}
\end{center}
where $n=(n_1,n_2)$ denotes the outer unit normal to $\partial \Omega$ and $\Gamma$ is an open portion of the boundary of $\Omega $. These boundary conditions occur when the constraint on $\Gamma$ prevents the transverse displacement of the plate and, simultaneously, leaves free the rotation of the transverse material fibers ${x} \times \left [ - \frac{h}{2}, \frac{h}{2}   \right ]$, $x \in \Gamma$, around the direction of the tangent to the boundary of the plate. Situations of this type are quite common in engineering structures, for example for slab plates or roof plates \cite{Timoshenko_Woinowsky-Krieger}. 

\noindent Assuming $B,\nu\in C^2(\overline{\Omega})$ and $ \Gamma$ of $C^{4,\alpha}$ class, we obtain the following  doubling inequality at the boundary for  problem \eqref{eq:equazione_piastra_intro}-- \eqref{eq:bdry_cond2_intro}
\begin{equation}
	\label{eq:doubling_intro}
	\int_{B_{s}(P)\cap \Omega}|u|^2\leq CN^k\left(s\over r \right)^{\log_2(CN^k)}\int_{B_{r}(P)\cap \Omega}|u|^2,
\end{equation}
for every $P\in\Gamma$ such that $B_{r_0}(P)\cap\partial\Omega\subset\Gamma$, where $C,k>1, 0<r<s< \frac{r_0}{C}$  and $N$ is the frequency of the solution $u$ defined as follows 
\begin{equation}
	\label{eq:N_intro}
	N=\frac{\int_{B_{r_0}(P)\cap \Omega}|u|^2}{\int_{B_{\frac{r_0}{C}}(P)\cap \Omega}|u|^2}.
\end{equation}

\noindent Inequalities as \eqref{eq:doubling_intro}  are classical tools of Strong Unique Continuation at the Boundary arising in the quantification of the local vanishing rate on the boundary of solutions to PDEs (see, for instance, \cite{AdE}, \cite{l:arv}, \cite{l:BaGa}, \cite{DCSV}, \cite{l:KuNy}, \cite{MRV-LeMatem-2020}). The quantification of this property poses significant challenges, in fact one needs to infer quantitative control on the zero set of $u$.  Moreover, from a more applied side, it has been proved useful in the context of stability estimates for inverse problems consisting in the determination of unknown coefficients and boundaries (\cite{ABRV, l:mrv19, MRV-JEMS, S}). 

\noindent While the overall strategy follows the ideas introduced in \cite{MRV-LeMatem-2020}, we emphasize that the presence of supported boundary conditions displays several new features and difficulties that we outline below.  As first step, following \cite{MRV-LeMatem-2020}, we flatten the boundary by means of a conformal map which has been introduced in \cite{l:arv}. Such a change of variables on one hand simplifies the geometry of the problem and preserves the structure of the equation, but on the other hand leads to a quite demanding computation in order to achieve the underlying boundary conditions. A further obstruction in dealing with supported boundary conditions is due to the fact that, in order to apply a reflection argument, we need to operate another transformation  (see \eqref{eq:19-2bis}--\eqref{eq:def di v} for a precise definition) that fits our problem into a new one for a solution $v$ satisfying more elementary boundary conditions, namely

\begin{equation*}
v=\Delta v=0, \quad \hbox{in } (-1,1)\times \{0\}.
\end{equation*}

\noindent The latter allows us to apply an odd reflection to $v$, with extended function $\overline{v}$ belonging to $H^4(B_1) $ and satisfying a suitable partial differential equation in $B_1$.

\noindent Taking advantage of the regularity of the extended solution $\overline{v}$, we are now in position to deal with the following Carleman estimate (\cite{MRV-LeMatem-2020}). For every $\tau\geq \overline{\tau}$ and for every $r\in (0, 1)$ the following inequality holds 

\begin{gather}
    \label{eq:catania24.4 intro}
	\tau^4r^2\int\rho^{-2-2\tau}|U|^2 +\sum_{k=0}^3 \tau^{6-2k}\int\rho^{2k+1-2\tau}|D^kU|^2 \\\nonumber  \leq C
	\int\rho^{8-2\tau}(\Delta^2 U)^2,
\end{gather}
 for every $U\in C^\infty_0(B_{1}\setminus\ \overline{B}_{r/4})$, where $\rho(x,y)\sim \sqrt{x^2+y^2}$ as $(x,y)\rightarrow (0,0)$ (see Proposition \ref{prop:Carleman} for a precise statement). Indeed we finally apply the  estimate \eqref{eq:catania24.4 intro} to $U=\eta\overline{v}$ where $\eta$ is a cut-off function and we conclude by using a standard procedure. We also remark that the Carleman estimate \eqref{eq:catania24.4 intro} for the bi-laplacian operator is analogous to the estimate $(6.32)$ derived in \cite{MRV-IUMJ-2007}, with the exception of the presence of the first extra integral on the left hand side of \eqref{eq:catania24.4 intro} (see \cite{Ba} for the first derivation of a doubling inequality in the interior {}from a Carleman estimate of the form \eqref{eq:catania24.4 intro}). It can be shown that this term plays a crucial role in obtaining the doubling inequality \eqref{eq:doubling_intro}.

\noindent The outline of the article is as follows. After briefly recalling our notation  in Section 2, in Section 3 we introduce and state our main result, Theorem \ref{theo:DoublingInequality}. In Section 4, we recall two main ingredients of our strategy, namely the  conformal mapping in Theorem \ref{prop:conf_map} and the Carleman estimate in Proposition \ref{prop:Carleman}. In Section 5, which constitutes the central estimates of the paper, we deduce the doubling inequality at the boundary. The proof is divided into three steps, namely the reduction to a flat boundary, the extension by odd reflection and the application of the Carleman estimate.

\section{Notation} 
\label{sec:notation}

We shall generally denote points in $\R^2$ by $x=(x_1,x_2)$ or $y=(y_1,y_2)$, except for
Sections \ref{sec:Preliminary} and \ref{sec:Proof Doubling} where we rename $x,y$ the coordinates in $\R^2$.

We shall denote by $B_r(P)$ the disc in $\R^2$ of radius $r$ and
center $P$, by $B_r$ the disk of radius $r$ and
center $O$, by $B_r^+$, $B_r^-$ the hemidiscs in $\R^2$  of radius $r$ and
center $O$ contained in the halfplanes $\R^2_+= \{x_2>0\}$, $\R^2_-= \{x_2<0\}$ respectively, and by $R_{
a,b}$ the rectangle $(-a,a)\times(-b,b)$.

Given a matrix $A =(a_{ij})$, we shall denote by $|A|$ its Frobenius norm $|A|=\sqrt{\sum_{i,j}a_{ij}^2}$.

\begin{definition}
  \label{def:reg_bordo} (${C}^{k,\alpha}$ regularity)
Let $\Omega$ be a bounded domain in ${\R}^{2}$. Given $k,\alpha$,
with $k\in\N$, $0<\alpha\leq 1$, we say that a portion $S$ of
$\partial \Omega$ is of \textit{class ${C}^{k,\alpha}$ with
constants $r_{0}$, $M_{0}>0$}, if, for any $P \in S$, there
exists a rigid transformation of coordinates under which we have
$P=0$ and
\begin{equation*}
  \Omega \cap R_{r_0,2M_0r_0}=\{x \in R_{r_0,2M_0r_0} \quad | \quad
x_{2}>g(x_1)
  \},
\end{equation*}
where $g$ is a ${C}^{k,\alpha}$ function on
$[-r_0,r_0]$
satisfying
\begin{equation*}
g(0)=g'(0)=0,
\end{equation*}
\begin{equation*}
\|g\|_{{C}^{k,\alpha}([-r_0,r_0])} \leq M_0r_0,
\end{equation*}
where
\begin{equation*}
\|g\|_{{C}^{k,\alpha}([-r_0,r_0])} = \sum_{i=0}^k  r_0^i\sup_{[-r_0,r_0]}|g^{(i)}|+r_0^{k+\alpha}|g|_{k,\alpha},
\end{equation*}
\begin{equation*}
|g|_{k,\alpha}= \sup_ {\overset{\scriptstyle t,s\in [-r_0,r_0]}{\scriptstyle
t\neq s}}\left\{\frac{|g^{(k)}(t) - g^{(k)}(s)|}{|t-s|^\alpha}\right\}.
\end{equation*}
\end{definition}
\noindent Throughout the paper, summation over repeated indexes is assumed.

\noindent The cartesian components of the divergence of a second order tensor field $T(x)$ are defined, as usual, by 
\begin{equation}
	\label{eq:def-div}
	({\rm {div}}\ T(x))_i= T_{ij,j}(x), \quad i=1,2.
\end{equation}

\noindent In places we will use equivalently the symbols $D$ and $\nabla$ to denote the gradient of a function. Moreover, setting $\beta=(\beta_1,\beta_2) \in \N^2$, the multi-index notation $D^\beta$ represents the partial derivative
\begin{equation}
	\label{eq:def-multiindex}
	\dfrac{\partial^{\beta_1 +\beta_2} }{\partial x_1^{\beta_1} \partial x_2^{\beta_2}}.
\end{equation}
Throughout the paper, $C$ shall denote a positive constant which may change {}from line to line, and $p_2$, $q_2$, $P_2$, $Q_2$ shall denote second order differential operators which may also change {}from line to line.

\section{Main results} 
\label{sec:main}

In this paper we consider an isotropic thin elastic plate $\Omega\times \left[-\frac{h}{2},\frac{h}{2}\right]$, having middle plane $\Omega$ and thickness $h$. According to the Kirchhoff-Love theory, the transverse displacement $u$ of the plate satisfies the following fourth-order partial differential equation
\begin{equation}
    \label{eq:equazione_piastra}
    L(u) := {\rm div}\left ({\rm div} \left ( B(1-\nu)\nabla^2 u + B\nu \Delta u I_2 \right ) \right )=0, \qquad\hbox{in
    } \Omega,
\end{equation}
where the \emph{bending stiffness} $B$ is given by
\begin{equation}
  \label{eq:3.stiffness}
  B(x)=\frac{h^3}{12}\left(\frac{E(x)}{1-\nu^2(x)}\right),
\end{equation}
with the \emph{Young's modulus} $E$ and the \emph{Poisson's coefficient} $\nu$  given by
\begin{equation}
  \label{eq:3.E_nu}
  E(x)=\frac{\mu(x)(2\mu(x)+3\lambda(x))}{\mu(x)+\lambda(x)},\qquad\nu(x)=\frac{\lambda(x)}{2(\mu(x)+\lambda(x))},
\end{equation}
where $\lambda$ and $\mu$ are the Lam\'{e} moduli. Therefore, we may also rewrite \eqref{eq:3.stiffness} as
\begin{equation}
	\label{eq:3.B_riscritto}
	B(x)=\frac{h^3}{3}\cdot\frac{\mu(x)(\mu(x)+\lambda(x))}{2\mu(x)+\lambda(x)}.
\end{equation}
We assume that the Lam\'{e} moduli satisfy the strong convexity assumptions 
\begin{equation}
  \label{eq:3.Lame_convex}
  \mu(x)\geq \alpha_0>0,\qquad 2\mu(x)+3\lambda(x)\geq\gamma_0>0, \qquad \hbox{ in } \Omega,
\end{equation}
where $\alpha_0$, $\gamma_0$ are positive constants. Let us notice that the above assumptions imply that
\begin{equation}
  \label{eq:3.mu+lambda}
  \mu(x)+\lambda(x)\geq \min\left\{\alpha_0, \frac{\gamma_0}{2}\right\}, \qquad \hbox{ in } \Omega.
\end{equation}
Moreover, we shall assume the following regularity on $\lambda$ and $\mu$:
\begin{equation}
	\label{eq:regularity-Lamè}
	\| \lambda \|_{ C^2( \overline{\Omega}_{ r_0   }     )     }, 
	\quad 
	\| \mu \|_{ C^2( \overline{\Omega}_{ r_0   }     )     } 
	\leq
	\Lambda_0,
\end{equation}
with $\Lambda_0$ a positive constant.

We found it convenient to represent the plate equation $L(u)$ in the compact form
\begin{equation}
    \label{eq:equazione_piastra_C}
    L(u) =\partial^2_{ij}(c_{ijlk}\partial^2_{lk}u)=0 \qquad\hbox{in
    } \Omega,
\end{equation}
with
\begin{equation}
    \label{eq:coefficienti_C}
    c_{ijlk}=B(1-\nu)\delta_{il}\delta_{jk}+B\nu\delta_{ij}\delta_{lk}, \qquad\hbox{in
    } \Omega,
\end{equation}
so that, by a straightforward computation, equation \eqref{eq:equazione_piastra} can be rewritten as
\begin{equation}
    \label{eq:equazione_piastra_con B}
    L(u) =B(\Delta^2 u + \widetilde{a} \cdot \nabla (\Delta u) + \widetilde{q}_2(u))=0, \quad \hbox{in } \Omega,
\end{equation}
where 
\begin{equation}
	\label{eq:vettore_a_tilde}
	\widetilde{a}=2\frac{\nabla B}{B},
\end{equation}
\begin{equation}
	\label{eq:q_2_tilde}
	\widetilde{q}_2(u)=\sum_{i,j=1}^2\frac{1}{B}\partial^2_{ij}(B(1-\nu)+\nu B\delta_{ij})\partial^2_{ij} u.
\end{equation}
Let us define
\begin{equation}
	\label{eq:Omega_r_0}
	\Omega_{r_0} = \Omega \cap R_{r_0,2M_0r_0}=\left\{ x\in R_{r_0,2M_0r_0}\ |\ x_2>g(x_1) \right\},
\end{equation}
\begin{equation}
	\label{eq:Gamma_r_0}
	\Gamma_{r_0} = \partial\Omega \cap R_{r_0,2M_0r_0}= \left\{(x_1,g(x_1))\ |\ x_1\in (-r_0,r_0)\right\},
\end{equation}
with $g$ as in Section \ref{sec:notation} and such that
\begin{equation}
	\label{eq:regolar-boundary}
	\|g\|_{C^{4,\alpha }( [ -r_0, r_0  ]  )   } \leq M_0 r_0,
\end{equation}
for some $\alpha \in (0,1)$. 

We notice that, by \eqref{eq:3.Lame_convex}--\eqref{eq:regularity-Lamè}, the coefficient $\widetilde{a}$ and the second order operator $\widetilde{q_2} =\sum_{|\alpha| \leq 2} c_\alpha D^\alpha$ satisfy
\begin{equation}
	\label{eq:regul-atilde-q2tilde}
	\| \widetilde{a}\|_{ C^1( \overline{\Omega}_{r_0}, \R^2   )   } \leq L, \quad	
	\| c_\alpha \|_{ C^0( \overline{\Omega}_{r_0}   )   } \leq L,
\end{equation}
with $L>0$ only depending on $\alpha_0$, $\gamma_0$, $\Lambda_0$.

We restrict our attention to the set $\Omega_{r_0}$, and we consider the following
boundary conditions on $\Gamma_{r_0}$:
\begin{center}
\( {\displaystyle \left\{
\begin{array}{lr}
     u=0, \quad &\hbox{on } \Gamma_{r_0},\\
c_{ijlk}\partial^2_{lk}u\,n_in_j=0, \quad &\hbox{on } \Gamma_{r_0},
\end{array}
\right. } \) \vskip -4.2em
\begin{eqnarray}
& & \label{eq:bdry_cond1}\\
& & \label{eq:bdry_cond2}
\end{eqnarray}
\end{center}
where $n=(n_1,n_2)$ denotes the outer unit normal to $\partial\Omega$. By setting
\begin{equation}
\label{eq:de_nn}
\partial^2_{nn}u = \partial^2_{ij}u\,n_in_j,
\end{equation}
the second boundary condition can be rewritten as
\begin{equation}
\label{eq:rewrite_second_bdry_cond}
c_{ijlk}\partial^2_{lk}u\,n_in_j =B(1-\nu)\partial^2_{nn}u +B\nu\Delta u=0, \quad 
\hbox{on } \Gamma_{r_0}.
\end{equation}
Therefore, we study the following problem
\begin{center}
\( {\displaystyle \left\{
\begin{array}{lr}
L(u)=0, \quad &\hbox{in } \Omega_{r_0},\\
     u=0, \quad &\hbox{on } \Gamma_{r_0},\\
B(1-\nu)\partial^2_{nn}u +B\nu\Delta u=0, \quad &\hbox{on } \Gamma_{r_0}.
\end{array}
\right. } \) \vskip -5.5em
\begin{eqnarray}
& & \label{eq:equation}\\
& & \label{eq:bdry_cond1bis}\\
& & \label{eq:bdry_cond2bis}
\end{eqnarray}
\end{center}
Let us derive the \emph{variational formulation} to the problem \eqref{eq:equation}--\eqref{eq:bdry_cond2bis}. Let
\begin{equation}
\label{eq:set_test}
{\mathcal H}=\{v\in C^\infty(\overline{\Omega}_{r_0})\ |\ D^\alpha v = 0 \hbox{ on} 
\ \partial\Omega_{r_0} \setminus \Gamma_{r_0}, \forall \alpha \in \N^2, v=0 
\hbox{ on}\ \Gamma_{r_0}\}
\end{equation}
A \emph{weak solution} to problem \eqref{eq:equation}--\eqref{eq:bdry_cond2bis} is a function $u\in H^2(\Omega_{r_0})$ satisfying
\begin{center}
\( {\displaystyle \left\{
\begin{array}{l}
     \int_{\Omega_{r_0}}c_{ijlk}\partial^2_{lk}u\partial^2_{ij}v=0, \quad \hbox{for every } v\in {\mathcal H},
\\
u=0, \quad \hbox{on } \Gamma_{r_0}.
\end{array}
\right. } \) \vskip -4.3em
\begin{eqnarray}
& & \label{eq:variational1}\\
& & \label{eq:variational2}
\end{eqnarray}
\end{center}
 We notice that the elliptic and regularity assumptions on the Lam\'{e} moduli \eqref{eq:3.Lame_convex} and \eqref{eq:regularity-Lamè}, and the regularity condition \eqref{eq:regolar-boundary} on $g$ guarantee that $u \in H^4(\Omega_{r_0} )$, see, for instance, \cite{Agmon-book}.  
 
Hence, by applying the integration by parts we have that, for every $v\in {\mathcal H}$,
\begin{multline*}
\int_{\Omega_{r_0}}c_{ijlk}\partial^2_{lk}u\partial^2_{ij}v =
\int_{\Omega_{r_0}}\partial_i(c_{ijlk}\partial^2_{lk}u\partial_{j}v)-
\partial_i(c_{ijlk}\partial^2_{lk}u)\partial_{j}v= \\
=\int_{\partial\Omega_{r_0}}c_{ijlk}\partial^2_{lk}u\partial_{j}v\,n_i-
\int_{\Omega_{r_0}}\partial_i(c_{ijlk}\partial^2_{lk}u)\partial_{j}v=\\
=\int_{\Gamma_{r_0}}c_{ijlk}\partial^2_{lk}u\partial_{j}v\,n_i+
\int_{\Omega_{r_0}}\partial^2_{ij}(c_{ijlk}\partial^2_{lk}u)v.
\end{multline*}      
Moreover, since $v\equiv 0$ on $\Gamma_{r_0}$, we have
\begin{equation*}
\nabla v=(\nabla v\cdot n)n, \quad \hbox{on } \Gamma_{r_0},
\end{equation*}
\begin{equation*}
\partial_j v=(\partial_n v)n_j, \ j=1,2, \quad \hbox{on } \Gamma_{r_0},
\end{equation*}
so that
\begin{equation}
\label{eq:variational3}
\int_{\Omega_{r_0}}c_{ijlk}\partial^2_{lk}u\partial^2_{ij}v 
=\int_{\Gamma_{r_0}}(c_{ijlk}\partial^2_{lk}u\,n_in_j)\partial_n v+
\int_{\Omega_{r_0}}\partial^2_{ij}(c_{ijlk}\partial^2_{lk}u)v.
\end{equation} 

Therefore, given a weak solution $u\in H^2(\Omega_{r_0})$ to problem \eqref{eq:equation}--\eqref{eq:bdry_cond2bis}, by choosing the test functions in $C^\infty_0(\Omega_{r_0})\subset {\mathcal H}$ in \eqref{eq:variational1}, we get the differential equation \eqref{eq:equazione_piastra_C} and, consequently, we obtain
\begin{equation*}
\int_{\Gamma_{r_0}}(c_{ijlk}\partial^2_{lk}u\,n_in_j)\partial_n v=0, \quad \forall v\in {\mathcal H},
\end{equation*} 
which implies the boundary condition \eqref{eq:bdry_cond2bis}. Viceversa, a classical solution to problem \eqref{eq:equation}--\eqref{eq:bdry_cond2bis}, in view of \eqref{eq:variational3}, must satisfy the weak formulation \eqref{eq:variational1}--\eqref{eq:variational2}.

We are now in position to state the main result of this paper.

\begin{theo}[Doubling inequality at the supported boundary]
\label{theo:DoublingInequality}
Under the above notation, let us assume that the Lam$\acute{e}$ moduli satisfy \eqref{eq:3.Lame_convex}--\eqref{eq:regularity-Lamè} and the boundary $\Gamma_{r_0}$ is of $C^{4,\alpha}$-class, with $g$ satisfying \eqref{eq:regolar-boundary} for some $\alpha \in (0,1)$. Then, there exist $k>1$ and $C>1$ only depending on $\alpha_0$, $\gamma_0$, $\Lambda_0$, $M_0$, $\alpha$, such that, for every $0<r<s<\frac{r_0}{C}$ we have that
\begin{equation}
	\label{eq:doubling}
	\int_{B_{s}\cap \Omega_{r_0}}|u|^2\leq CN^k\left(s\over r \right)^{\log_2(CN^k)}\int_{B_{r}\cap \Omega_{r_0}}|u|^2,
\end{equation}
where
\begin{equation}
	\label{eq:N}
	N=\frac{\int_{B_{r_0}\cap \Omega_{r_0}}|u|^2}{\int_{B_{\frac{r_0}{C}}\cap \Omega_{r_0}}|u|^2}.
\end{equation}
\end{theo}

\section{Preliminary results} 
\label{sec:Preliminary}

The following Proposition, obtained in \cite{l:arv}, introduces a conformal map which flattens the boundary $\Gamma_{r_0}$ and preserves the structure of equation \eqref{eq:equazione_piastra_con B}.

\begin{prop} [{\bf Conformal mapping} (\cite{l:arv}, Proposition 3.1)]
    \label{prop:conf_map}
Under the hypotheses of Theorem \ref{theo:DoublingInequality}, there exists an injective sense preserving differentiable map
\begin{equation*}
	\Phi=(\varphi,\psi):[-1,1]\times[0,1]\rightarrow
	\overline{\Omega}_{r_0}
\end{equation*}
which is conformal and satisfies
\begin{equation}
	\label{eq:9.assente}
	\Phi((-1,1)\times(0,1))\supset B_{\frac{r_0}{K}}(0)\cap \Omega_{r_0},
\end{equation}
\begin{equation}
	\label{eq:9.2b}
	\Phi(([-1,1]\times\{0\})= \left\{ (x_1,g(x_1))\ |\ x_1\in [-r_1,r_1]\right\} \subset \Gamma_{r_0},
\end{equation}

\begin{equation}
	\label{eq:9.2a}
	\Phi(0,0)= (0,0),
\end{equation}

\begin{equation}
	\label{eq:gradPhi}
	\frac{c_0r_0}{2C_0}\leq |D\Phi(y)|\leq \frac{r_0}{2}, \quad \forall y\in [-1,1]\times[0,1],
\end{equation}
\begin{equation}
	\label{eq:gradPhiInv}
	\frac{4}{r_0}\leq |D\Phi^{-1}(x)|\leq \frac{4C_0}{c_0r_0},
	\quad\forall x\in \Phi([-1,1]\times[0,1]),
\end{equation}
\begin{equation}
	\label{eq:stimaPhi}
	\frac{r_0}{K}|y|\leq|\Phi(y)|\leq \frac{r_0}{2}|y|, \quad \forall
	y\in [-1,1]\times[0,1],
\end{equation}
with $K>8$, $0<c_0<C_0$, $\frac{r_1}{r_0}$ being constants only depending on $M_0$
and $\alpha$.
\end{prop}

Another basic ingredient for our proof of the doubling inequality at the supported boundary of the plate is the following Carleman estimate, whose proof can be found in \cite[Proposition $3.5$]{MRV-LeMatem-2020}.

\begin{prop} [{\bf Carleman estimate}]
    \label{prop:Carleman}
		Let us define
\begin{equation}
    \label{eq:catania24.2}
		\rho(x,y) = \phi\left(\sqrt{x^2+y^2}\right),
\end{equation}
where
\begin{equation}
    \label{eq:catania24.3}
		\phi(s) = \frac{s}{\left(1+\sqrt{s}\right)^2}.
\end{equation}
Then there exist absolute constants $\overline{\tau}>1$, $C>1$ such that
\begin{gather}
    \label{eq:catania24.4}
	\tau^4r^2\int\rho^{-2-2\tau}|U|^2 +\sum_{k=0}^3 \tau^{6-2k}\int\rho^{2k+1-2\tau}|D^kU|^2 \\\nonumber  \leq C
	\int\rho^{8-2\tau}(\Delta^2 U)^2,
\end{gather}
for every $\tau\geq \overline{\tau}$, for every $r\in (0, 1)$ and for every $U\in C^\infty_0(B_{1}\setminus\ \overline{B}_{r/4})$.
\end{prop}

\section{Proof of the doubling inequality at the boundary} 
\label{sec:Proof Doubling}

The strategy of the proof consists of several steps. Firstly, we flatten the boundary $\Gamma_{r_0}$ by means of the conformal mapping defined in Proposition \ref{prop:conf_map}. The conformal mapping preserves the structure of the operator  and, by the peculiar supported boundary conditions, the solution can be extended by an odd reflection with respect to the flattened boundary. Then, the desired doubling inequality \eqref{eq:doubling} follows {}from the Carleman estimate \eqref{eq:catania24.4} and a standard argument.

{}From now on, we shall simply denote by $R$ the rectangle $(-1,1)\times(0,1)$.

\medskip

\textit{First step. Reduction to a flat boundary.}

\medskip

Given a weak solution $u\in H^2(\Omega_{r_0})$ to problem \eqref{eq:equation}--\eqref{eq:bdry_cond2bis}, let us denote

\begin{equation*}
w(y)=u(\Phi(y)),
\end{equation*}
where $\Phi=(\Phi_1,\Phi_2)=(\varphi, \psi)$ is the conformal mapping introduced in Proposition \ref{prop:conf_map}.

It is useful to notice that, since $\Phi$ is a conformal map,
\begin{equation*}
   \label{eq:DPhi}
  D\Phi 
	=\left( \begin{array}{ll}
  \partial_{y_1}\Phi_1 &\partial_{y_2}\Phi_1\\
   &  \\
  \partial_{y_1}\Phi_2 &\partial_{y_2}\Phi_2\\
  \end{array}\right)
	=\left( \begin{array}{ll}
  \partial_{y_1}\varphi & \partial_{y_2}\varphi\\
   &  \\
  -\partial_{y_2}\varphi &\partial_{y_1}\varphi\\
  \end{array}\right),
\end{equation*}
\begin{equation*}
  det(D\Phi(y)) = |\nabla\varphi(y)|^2,
\end{equation*}
\begin{equation*}
  (D\Phi)^{-1} =\frac{1}{|\nabla \varphi|^2}\left( \begin{array}{ll}
  \varphi_{y_1} &-\varphi_{y_2}\\
   &  \\
  \varphi_{y_2} &\varphi_{y_1}\\
  \end{array}\right).
\end{equation*}
For any function $z=z(x)$ defined in $\Phi(R)$, we can compute
\begin{equation}
	\label{eq:6.29 di arv}
   (\nabla_x z) (\Phi(y)) = [(D\Phi(y))^{-1}]^T\nabla_y (z(\Phi(y)).
\end{equation}
It follows that
\begin{equation}
		\label{eq:6.32 di arv}
   (\Delta u) (\Phi(y)) = \frac{1}{|\nabla \varphi(y)|^2}(\Delta w)(y)
\end{equation}
and
\begin{multline}
		\label{eq:6.34 di arv}
   (\Delta^2 u) (\Phi(y)) = \frac{1}{|\nabla \varphi(y)|^2}\Delta \left(\frac{1}{|\nabla \varphi(y)|^2}\Delta w(y)\right)=\\
	=\frac{1}{|\nabla \varphi(y)|^2}\left[\frac{1}{|\nabla \varphi(y)|^2}\Delta^2 w
	+2\nabla\left(\frac{1}{|\nabla \varphi(y)|^2}\right)\cdot\nabla(\Delta w)+q_2(w)\right],
\end{multline}
where $q_2$ is a second order differential operator with $C^1$ norm of its coefficients bounded in terms of $\alpha_0$, $\gamma_0$, $\Lambda_0$.

In order to derive the differential equation satisfied by $w$ in $R$, let us consider for the moment test functions having compact support in $R$. Precisely, given $\widetilde{V}\in C^\infty_0(R)$, let
\begin{equation*}
\widetilde{v}(x)=\widetilde{V}(\Phi^{-1}(x)), \qquad \hbox {in  } \Phi(R)\subset \Omega_{r_0}.
\end{equation*}
By integrating by parts \eqref{eq:variational1} and recalling \eqref{eq:equazione_piastra_con B}, we have
\begin{multline*}
0=\int_{\Omega_{r_0}} c_{ijlk}\partial^2_{lk}u\partial^2_{ij}\widetilde{v}=\int_{\Phi(R)} c_{ijlk}\partial^2_{lk}u\partial^2_{ij}\widetilde{v}=
\int_{\Phi(R)} \partial^2_{ij}(c_{ijlk}\partial^2_{lk}u)\widetilde{v}=\\
= \int_{\Phi(R)} B(\Delta^2 u + \widetilde{a} \cdot \nabla (\Delta u) + \widetilde{q}_2(u))
\widetilde{v},
\end{multline*}
where $\widetilde{a}$ and $\widetilde{q}_2(u)$ are defined in \eqref{eq:vettore_a_tilde} and  \eqref{eq:q_2_tilde}, respectively.

By operating the change of variables $x=\Phi(y)$ in the last integral, we have that, for every $\widetilde{V}\in C^\infty_0(R)$,
\begin{multline*}
\int_{R} B(\Phi(y))\left [ (\Delta^2 u)(\Phi(y)) + \widetilde{a}(\Phi(y))\cdot (\nabla (\Delta u))(\Phi(y)) + ( \widetilde{q}_2  (u))(\Phi(y))\right]\cdot\\
|det(D\Phi(y)|\widetilde{V}(y)=0.
\end{multline*}
Since $det (D\Phi(y))=|\nabla \varphi(y)|^2$, by the arbitrariness of the test functions $\widetilde{V}\in C^\infty_0(R)$, we obtain that
\begin{equation}
   \label{eq:29-1}
(\Delta^2 u)(\Phi(y)) + \widetilde{a}(\Phi(y))\cdot (\nabla (\Delta u))(\Phi(y)) + ( \widetilde{q}_2(u))(\Phi(y))=0, 
\qquad \hbox{ in } R. 
\end{equation}
%
%
%
By applying \eqref{eq:6.29 di arv} with $z=\Delta u$, and by \eqref{eq:6.32 di arv}, we compute
\begin{multline}
   \label{eq:30-3}
(\nabla_x (\Delta u)) (\Phi(y)) = [(D\Phi(y))^{-1}]^T\nabla_y ((\Delta u)(\Phi(y))=\\
=[(D\Phi(y))^{-1}]^T\nabla_y \left(\frac{1}{|\nabla \varphi|^2}\Delta w(y)\right)
=[(D\Phi(y))^{-1}]^T\frac{1}{|\nabla \varphi|^2}\nabla_y \Delta w(y)+q_2(w),
\end{multline}
where $q_2$ is a second order differential operator with $C^2$-norm of its coefficients only depending on $M_0, \alpha, \alpha_0, \gamma_0, \Lambda_0$.

By \eqref{eq:29-1}, \eqref{eq:6.32 di arv} and \eqref{eq:30-3}, we have 
\begin{multline}\label{eq:30-4}
\frac{1}{|\nabla \varphi|^2} \Delta^2 w + 
2\nabla\left(\frac{1}{|\nabla \varphi|^2}\right)\cdot \nabla_y (\Delta w) +  \widetilde{a}(\Phi(y))\cdot [(D\Phi(y))^{-1}]^T\nabla_y (\Delta w) + q_2(w) =
\\
\frac{1}{|\nabla \varphi|^2} \Delta^2 w + 
\left(2\nabla\left(\frac{1}{|\nabla \varphi|^2}\right)+
 [(D\Phi(y))^{-1}] \widetilde{a}(\Phi(y))\right)
\cdot \nabla_y (\Delta w) + q_2(w)=0, 
\\
\qquad \hbox{in } R,
\end{multline}
where $q_2$ is a second order differential operator with $C^0$-norm on $\overline{R}$ of its coefficients only depending on $M_0, \alpha, \alpha_0, \gamma_0, \Lambda_0$. Therefore the function $w$ satisfies
\begin{equation*}
{\mathcal L}(w)=0, \qquad \hbox{ in } R,
\end{equation*}
with
\begin{equation}
\label{eq:27-1}
	{\mathcal L}(w)=\Delta^2 w+b\cdot \nabla (\Delta w)+Q_2(w),
\end{equation}
where $b \in C^1 ( \overline{R}, \R^2)$ is a vector field and $Q_2= \sum_{|\alpha| \leq 2} C_\alpha D^\alpha$ is a second order differential operator such that $\| b\|_{  C^1 ( \overline{R}, \R^2)  } $ and $\| C_\alpha\|_{  C^0 ( \overline{R})} $ are bounded in terms of $M_0, \alpha, \alpha_0, \gamma_0, \Lambda_0$.

In order to derive the boundary conditions satisfied by $w=u\circ\Phi$, with $u\in H^2(\Omega_{r_0})$ a weak solution to problem \eqref{eq:equation}--\eqref{eq:bdry_cond2bis}, let us define
\begin{equation}
\label{eq:set_test_R}
\widetilde{{\mathcal H}}=\{v\in C^\infty(\Phi(R))\ |\ D^\alpha v = 0 \hbox{ on} 
\ \partial(\Phi(R)) \setminus (\Gamma_{r_0}\cap \Phi(R)), \forall \alpha \in \N^2, v=0 
\hbox{ on}\ \Gamma_{r_0}\cap \Phi(R)\}
\end{equation}
It follows that  
\begin{equation}
	\label{eq:7-4}
     \int_{\Phi(R)}c_{ijlk}\partial^2_{lk}u\partial^2_{ij}v=0,\forall v\in \widetilde{{\mathcal H}}.
\end{equation}

Given
$v\in \widetilde{{\mathcal H}}$, let us denote

\begin{equation*}
V(y)=v(\Phi(y)).
\end{equation*}

By operating the change of variables $x=\Phi(y)$, we have that, for every 
$v\in \widetilde{{\mathcal H}}$,

\begin{equation}
	\label{eq:7-4bis}
     0=\int_{\Phi(R)}c_{ijlk}\partial^2_{lk}u\partial^2_{ij}v=
		\int_{R}c_{ijlk}(\Phi(y))|det(D\Phi(y))|(\partial^2_{lk}u)(\Phi(y))(\partial^2_{ij}v)(\Phi(y)).
\end{equation}

Let us set
\begin{equation}
	\label{eq:8-2}
     M(y)=(D\Phi(y))^{-1}
\end{equation}

By applying \eqref{eq:6.29 di arv} with $z=\partial_{x_k}u$, we obtain

\begin{multline*}
     (\partial^2_{x_lx_k}u)(\Phi(y))=\{M^T(y)\nabla_y[(\partial_{x_k}u)(\Phi(y))]\}_l=\\
		=\{M^T(y)\nabla_y[M^T(y)\nabla_y w]_k\}_l= m_{hl}(y)\partial_{y_h}(m_{sk}(y)\partial_{y_s}w), \qquad l,k=1,2.
\end{multline*}

Similarly, given
$v\in \widetilde{{\mathcal H}}$, we have

\begin{equation*}
     (\partial^2_{x_ix_j}v)(\Phi(y))= m_{h'i}(y)\partial_{y_{h'}}(m_{s'j}(y)\partial_{y_{s'}}V), \qquad l,k=1,2.
\end{equation*}

{}From \eqref{eq:7-4bis}, for every $v\in \widetilde{{\mathcal H}}$, we have

\begin{multline}
	\label{eq:9-1}
     0=\int_{\Phi(R)}c_{ijlk}\partial^2_{lk}u\partial^2_{ij}v=\\
		=\int_{R}c_{ijlk}(\Phi(y))|det(D\Phi(y))|
		m_{hl}\partial_{y_h}(m_{sk}\partial_{y_s}w)
		m_{h'i}\partial_{y_{h'}}(m_{s'j}\partial_{y_{s'}}V)=\\
		=\int_{R}\widetilde{c}_{h'jhk}
		\partial_{y_h}(m_{sk}\partial_{y_s}w)
		\partial_{y_{h'}}(m_{s'j}\partial_{y_{s'}}V),
\end{multline}
where 
\begin{equation}
	\label{eq:9-2}
     \widetilde{c}_{h'jhk}(y)=m_{h'i}(y)c_{ijlk}(\Phi(y))|det(D\Phi(y))|m_{hl}(y).
\end{equation}

By twice integrating by parts, we obtain
\begin{equation}
	\label{eq:10-0}
     \int_{\partial R} {\mathcal B}_1- {\mathcal B}_2
		+\int_R\partial_{y_{s'}}[m_{s'j}\partial_{y_{h'}}(\widetilde{c}_{h'jhk}\partial_{y_h}(m_{sk}\partial_{y_s}w))]V=0,
\end{equation}
where
\begin{equation}
	\label{eq:10-1a}
     {\mathcal B}_1= n_{h'}\widetilde{c}_{h'jhk}\partial_{y_h}(m_{sk}\partial_{y_s}w)
		m_{s'j}\partial_{y_{s'}}V , 
\end{equation}

\begin{equation}
	\label{eq:10-1b}
     {\mathcal B}_2= n_{s'}m_{s'j}\partial_{y_{h'}}
		(\widetilde{c}_{h'jhk}\partial_{y_h}(m_{sk}\partial_{y_s}w))V.
\end{equation}

Next, by the definition of $\widetilde{{\mathcal H}}$, ${\mathcal B}_2\equiv 0$ on $\partial R$ and 

\begin{equation}
	\label{eq:10-2}
     \int_{\partial R} {\mathcal B}_1= \int_{[-1,1]\times\{0\}} {\mathcal B}_1=
		\int_{[-1,1]\times\{0\}} -(\widetilde{c}_{2jhk}\partial_{y_h}(m_{sk}\partial_{y_s}w)m_{2j})\partial_{y_2}V\ . 
\end{equation}

Let us compute, by using \eqref{eq:coefficienti_C},
\begin{multline}
	\label{eq:11-*}
	\widetilde{c}_{2jhk}\partial_{y_h}(m_{sk}\partial_{y_s}w)m_{2j}=
	|det(D\Phi(y))|m_{2i}m_{2j}m_{hl}c_{ijlk}
		\partial_{y_h}(m_{sk}\partial_{y_s}w)=\\
		=|det(D\Phi(y))|m_{2i}m_{2j}m_{hl}[B(1-\nu)\delta_{il}\delta_{jk}+B\nu\delta_{ij}\delta_{lk}]\partial_{y_h}(m_{sk}\partial_{y_s}w)=\\
		=B|det(D\Phi(y))|[m_{2l}m_{2k}m_{hl}(1-\nu)+m_{2i}m_{2i}m_{hk}\nu]
		(m_{sk}\partial^2_{y_sy_h}w+\partial_{y_h}(m_{sk})\partial_{y_s}w)=\\
		=
		B|det(D\Phi(y))|({\mathcal A}_1+{\mathcal A}_2),
\end{multline}
where
\begin{equation}
	\label{eq:A_1}
	{\mathcal A}_1=[(1-\nu)m_{2l}m_{2k}m_{hl}m_{sk}+\nu m_{2i}m_{2i}m_{hk}m_{sk}]\partial^2_{y_sy_h}w,
\end{equation}
\begin{equation}
	\label{eq:A_2}
	{\mathcal A}_2=[(1-\nu)m_{2l}m_{2k}m_{hl}\partial_{y_h}(m_{sk})+\nu m_{2i}m_{2i}m_{hk}\partial_{y_h}(m_{sk})]\partial_{y_s}w.
\end{equation}

In order to evaluate ${\mathcal A}_1$, let us notice that

\begin{equation}
	\label{eq:12-2}
	MM^T=\frac{1}{|\nabla \varphi|^2}I_2.
\end{equation}
Therefore,
\begin{equation}
	\label{eq:13-1}
	m_{2l}m_{2k}m_{hl}m_{sk}=\frac{1}{|\nabla \varphi|^4}\delta_{2h}\delta_{2s},
\end{equation}
\begin{equation}
	\label{eq:13-2}
	m_{2i}m_{2i}m_{hk}m_{sk}=\frac{1}{|\nabla \varphi|^4}\delta_{hs}.
\end{equation}
By inserting \eqref{eq:13-1} and \eqref{eq:13-2} in \eqref{eq:A_1}, we have
\begin{equation}
	\label{eq:13-3}
	{\mathcal A}_1=\frac{1}{|\nabla \varphi|^4}
	[(1-\nu)\delta_{2h}\delta_{2s}+\nu \delta_{hs}]\partial^2_{y_sy_h}w=
	\frac{1}{|\nabla \varphi|^4}
	[(1-\nu)\partial^2_{y_2 y_2}w+\nu \Delta w].
\end{equation}
On the other hand, since $w(y_1,0)\equiv 0$, we may rewrite \eqref{eq:13-3} as
\begin{equation}
	\label{eq:14-1}
	{\mathcal A}_1=
	\frac{1}{|\nabla \varphi|^4}
	\Delta w.
\end{equation}
Again by $w(y_1,0)\equiv 0$, we compute
\begin{multline}
	\label{eq:14-2}
	m_{2l}m_{2k}m_{hl}\partial_{y_h}(m_{sk})\partial_{y_s}w=(m_{2l}m_{hl})m_{2k}
	\partial_{y_h}(m_{2k})\partial_{y_2}w=\\
	=\frac{1}{2}\frac{1}{|\nabla \varphi|^2}\delta_{2h}\partial_{y_h}(m_{2k}m_{2k})\partial_{y_2}w=\frac{1}{2}\frac{1}{|\nabla \varphi|^2}\partial_{y_2}\left(\frac{1}{|\nabla \varphi|^2}\right)\partial_{y_2}w .
\end{multline}
Next, let us notice that
\begin{equation}
	\label{eq:15-2}
	(D\Phi)^T(D\Phi)=|\nabla \varphi|^2 I_2.
\end{equation}
By using again $w(y_1,0)\equiv 0$, by recalling that $\Phi_k$ are harmonic functions, $k=1,2$, and by \eqref{eq:15-2}, we compute
\begin{multline}
	\label{eq:16-1}
	m_{2i}m_{2i}m_{hk}\partial_{y_h}(m_{sk})\partial_{y_s}w
	=\frac{1}{|\nabla \varphi|^2}m_{hk}\partial_{y_h}(m_{2k})\partial_{y_2}w=\\=
	\frac{1}{|\nabla \varphi|^4}\partial_{y_h}(\Phi_k)\partial_{y_h}\left(\frac{\partial_{y_2}(\Phi_k)}{|\nabla \varphi|^2}\right)\partial_{y_2}w=\\
	=\frac{1}{|\nabla \varphi|^4}
	\left[\partial_{y_h}\left(\frac{\partial_{y_h}(\Phi_k)\partial_{y_2}(\Phi_k)}{|\nabla \varphi|^2}\right)-\frac{\Delta \Phi_k\partial_{y_2}(\Phi_k)}{|\nabla \varphi|^2}\right]\partial_{y_2}w=\\
	\frac{1}{|\nabla \varphi|^4}
	\left[\partial_{y_h}\left(\frac{\partial_{y_h}(\Phi_k)\partial_{y_2}(\Phi_k)}{|\nabla \varphi|^2}\right)\right]\partial_{y_2}w=
	\frac{1}{|\nabla \varphi|^4}
	\partial_{y_h}(\delta_{h2})\partial_{y_2}w=0.
\end{multline}
{}From \eqref{eq:14-2} and \eqref{eq:16-1}, we have
\begin{equation}
	\label{eq:16-2}
	{\mathcal A}_2=\frac{1}{2}(1-\nu)
	\frac{1}{|\nabla \varphi|^2}
	\partial_{y_2}\left(\frac{1}{|\nabla \varphi|^2}\right)
	\partial_{y_2}w.
\end{equation}
By \eqref{eq:10-2}, \eqref{eq:14-1} and \eqref{eq:16-2}, we can compute ${\mathcal B}_1$
on $[-1,1]\times \{0\}$:
\begin{multline}
	\label{eq:16-3}
{\mathcal B}_1=-(\widetilde{c}_{2jhk}\partial_{y_h}(m_{sk}\partial_{y_s}w)m_{2j})\partial_{y_2}V=\\
=-\frac{B}{|\nabla \varphi|^2}\left[\Delta w +\frac{1-\nu}{2}|\nabla \varphi|^2\partial_{y_2}\left(\frac{1}{|\nabla \varphi|^2}\right)\partial_{y_2}w\right]\partial_{y_2}V.
\end{multline}

By \eqref{eq:10-0}, we have that 
\begin{multline}
	\label{eq:17-1}
	\int_R\partial_{y_{s'}}[m_{s'j}\partial_{y_{h'}}(\widetilde{c}_{h'jhk}\partial_{y_h}(m_{sk}\partial_{y_s}w))]V=\\
	=\int_{[-1,1]\times \{0\}}\frac{B}{|\nabla \varphi|^2}\left[\Delta w +\frac{1-\nu}{2}|\nabla \varphi|^2\partial_{y_2}\left(\frac{1}{|\nabla \varphi|^2}\right)\partial_{y_2}w\right]\partial_{y_2}V.
\end{multline}
for every $V=v(\Phi(y))$ with $v\in \widetilde{{\mathcal H}}$ . 
By choosing the test functions $v\in C^\infty_0(\Phi(R))$, that is $V\in C^\infty_0(R)$, we obtain that $w$ satisfies the differential equation 
\begin{equation*}
\partial_{y_{s'}}[m_{s'j}\partial_{y_{h'}}(\widetilde{c}_{h'jhk}\partial_{y_h}(m_{sk}\partial_{y_s}w))]=0,
\end{equation*}
which coincides with the equilibrium equation \eqref{eq:30-4}, and therefore 
\begin{equation*}
\int_{[-1,1]\times \{0\}}\frac{B}{|\nabla \varphi|^2}\left[\Delta w +\frac{1-\nu}{2}|\nabla \varphi|^2\partial_{y_2}\left(\frac{1}{|\nabla \varphi|^2}\right)\partial_{y_2}w\right]\partial_{y_2}V=0 \ , 
\end{equation*}
for every $V\in C^\infty_0(R)$. By recalling \eqref{eq:27-1} and the boundary condition implied by the last identity, we have that $w=u(\Phi(y)) \in H^4(R)$ satisfies the following problem
\begin{center}
\( {\displaystyle \left\{
\begin{array}{lr}
\Delta^2 w+b\cdot\nabla\Delta w+Q_2(w)=0, \quad &\hbox{in } R,\\
     w(y_1,0)=0, \quad &\hbox{in } [-1,1]\times\{0\},\\
\Delta w+\gamma(y_1)\partial_{y_2}w=0, \quad &\hbox{in } [-1,1]\times\{0\},
\end{array}
\right. } \) \vskip -5.0em
\begin{eqnarray}
& & \label{eq:19-2a}\\
& & \label{eq:19-2b}\\
& & \label{eq:19-2c}
\end{eqnarray}
\end{center}
where
\begin{equation}
	\label{eq:19-2bis}
	\gamma(y_1)=\frac{1-\nu}{2}|\nabla \varphi|^2\partial_{y_2}\left(\frac{1}{|\nabla \varphi|^2}\right) \big |_{(y_1,0)}.
\end{equation}
In order to get a simpler boundary condition, let us introduce the function
\begin{equation}
    \label{eq:def di v}
v=e^{\frac{1}{2}y_2\gamma(y_1)}w.
\end{equation}
Denoting 
\begin{equation}
	\label{eq:a_def}
	a=e^{-\frac{1}{2}y_2\gamma(y_1)} , 
\end{equation}
we have
\begin{equation}
	\label{eq:w_v}
	w=av,
\end{equation}
\begin{equation}
	\label{eq:bilapl_w}
	\Delta^2w=\Delta^2(av) = a\Delta^2 v+4\nabla a\cdot \nabla(\Delta v) + p_2(v),
\end{equation}
\begin{equation}
	\label{eq:grad_lapl_w}
	\nabla(\Delta w)=a\nabla(\Delta v)+ p_2(v).
\end{equation}
Substituting \eqref{eq:w_v}, \eqref{eq:bilapl_w}, \eqref{eq:grad_lapl_w} in equation
\eqref{eq:27-1}, we obtain
\begin{equation}
	\label{eq:equa_v}
	\Delta^2 v+(\frac{4}{a}\nabla a + b)\cdot \nabla(\Delta v)+ p_2(v)=0.
\end{equation}
Noticing that $v(y_1,0)= 0$, $\partial_{y_1}v(y_1,0)=0$ and $a(y_1, 0)=1$ for every $y_1\in [-1,1]$, we can compute for every $(y_1,y_2)\in [-1,1]\times \{0\}$,
\begin{multline}
	\label{eq:bordo_v}
	\Delta w+\gamma(y_1)\partial_{y_2} w= \Delta(av)+\gamma(y_1)\partial_{y_2} v=
	2\nabla a \cdot \nabla v +\Delta v +\gamma(y_1)\partial_{y_2} v=\\
	=2 \partial_{y_2} a\partial_{y_2} v  +\Delta v +\gamma(y_1)\partial_{y_2} v=
	\Delta v.
\end{multline}
Therefore $v$ satisfies the following problem
\begin{center}
\( {\displaystyle \left\{
\begin{array}{lr}
\Delta^2 v+(\frac{4}{a}\nabla a + b)\cdot \nabla(\Delta v)+ p_2(v)=0, \quad &\hbox{in } R,\\
     v(y_1,0)=0, \quad &\hbox{in } [-1,1]\times\{0\},\\
\Delta v(y_1,0)=0, \quad &\hbox{in } [-1,1]\times\{0\} , 
\end{array}
\right. } \) \vskip -5.0em
\begin{eqnarray}
& & \label{eq:21-2a}\\
& & \label{eq:21-2b}\\
& & \label{eq:21-2c}
\end{eqnarray}
\end{center}
with $C^1$ and $C^0$ norm on $R$ of $a$ and of the coefficients of $p_2$ respectively depending on $M_0, \alpha, \alpha_0, \gamma_0, \Lambda_0$ only. 
\medskip

\textit{Second step. Extension by odd reflection.}

\medskip

Let us introduce the following extension of $v$ to $B_1$:
\begin{center}
	\( {\displaystyle \overline{v}(x,y) = \left\{
		\begin{array}{lr}
			v(x,y), \quad &\hbox{in } B_1^+,\\
			-v(x,-y), \quad &\hbox{in } B_1^-.
		\end{array}
		\right. } \) \vskip -3.0em
	\begin{eqnarray}
		& & \label{eq:ext-v}
	\end{eqnarray}
\end{center}
Moreover, let 
\begin{equation}
	\label{eq:def-f}
	f = - \left ( \frac{1}{a}(4\nabla a + ab)\cdot \nabla(\Delta v)+ p_2(v) \right ), \quad f \in L^2(B_1^+),
\end{equation}
and let us define
\begin{center}
	\( {\displaystyle \overline{f}(x,y) = \left\{
		\begin{array}{lr}
			f(x,y), \quad &\hbox{in } B_1^+,\\
			-f(x,-y), \quad &\hbox{in } B_1^-.
		\end{array}
		\right. } \) \vskip -3.0em
	\begin{eqnarray}
		& & \label{eq:ext-f}
	\end{eqnarray}
\end{center}
Then $\overline{f} \in L^2(B_1)$, $\overline{v} \in H^4(B_1)$ and 
\begin{equation}
	\label{eq:equaz-nel-disco}
	\Delta^2 \overline{v} = \overline{f}, \quad \hbox{in } B_1.
\end{equation}
The proof of \eqref{eq:equaz-nel-disco} can be obtained by adapting the arguments used in the proof of Proposition $4.1$ in \cite{l:arv}.

\medskip

\textit{Third step. Application of Carleman's estimate and conclusion.}

\medskip

Next, by a density argument, we apply the Carleman estimate \eqref{prop:Carleman} to $U=\eta(\sqrt{x^2+y^2}) \overline{v}$, where $\eta \in C_0^\infty ((0,1))$ is a suitable cut-off function; see, for instance, the proof of Lemma $4.1$ in \cite{MRV-LeMatem-2020}. We obtain the following result: there exists a positive number $\overline{r}_0 \in (0,1)$, only depending on $M_0, \alpha, \alpha_0, \gamma_0, \Lambda_0$, such that, for every $\overline{r}$ and for every $r$ such that $0 < 2r < \overline{r} < \dfrac{\overline{r}_0}{2}$, we have
\begin{multline}
	\label{eq:quasi-doubling}
	\overline{r}(2r)^{-2\tau} \int_{B_{2r}^+} |v|^2 + \overline{r}^{1-2\tau} \int_{B_{\overline{r}}^+} |v|^2 \leq\\
	\leq
	C
	\left (
	\left ( \dfrac{r}{4}    \right )^{-2\tau}  \int_{B_{r}^+} |v|^2
	+
	\left ( \dfrac{\overline{r}_0}{2}    \right )^{-2\tau}
	\int_{B_{\overline{r}_0}^+} |v|^2
	\right ),
\end{multline}
for every $\tau \geq \overline{\tau}$, with $\overline{\tau}$ a positive absolute constant and $C$ a positive constant depending on $M_0, \alpha, \alpha_0, \gamma_0, \Lambda_0$ only. 

Finally, by using a standard procedure (see, for instance, the proof of Theorem $2.2$ in \cite{MRV-LeMatem-2020}) and coming back to the original coordinates via Proposition \ref{prop:conf_map}, we obtain the desired doubling inequality  \eqref{eq:doubling} at the supported boundary of the plate.


\begin{thebibliography}{} \label{bbiibb}



\bibitem [AdE]{AdE}V. Adolfsson, L. Escauriaza, $C^{1,\alpha}$ domains and unique continuation at the boundary, Comm. Pure Appl. Math. L (1997), 935--969.

\bibitem[Ag]{Agmon-book}
S. Agmon, \emph{Lectures on Elliptic Boundary Value Problems}, Van Nostrand, New York, 1965.

\bibitem [ABRV]{ABRV} G. Alessandrini, E. Beretta, E. Rosset, S. Vessella, Optimal stability for inverse elliptic boundary value problems with unknown boundaries, Ann. Scuola Norm. Sup. Pisa Cl. Sci. (4) XXIX (2000), 755--806.

\bibitem[ARV]{l:arv} 
G. Alessandrini, E. Rosset, S. Vessella, Optimal three spheres inequality at the
boundary for the Kirchhoff-Love plate's equations with Dirichlet conditions, Arch. Ration. Mech. Anal. 231 (2019), 1455--1486.

\bibitem[Ba]{Ba} L. Bakri, Quantitative uniqueness for Schr\"{o}dinger operator, Indiana Univ. Math. J. 61 (4) (2012), 1565--1580.

\bibitem[BG]{l:BaGa} A. Banerjee, N. Garofalo, Quantitative uniqueness for elliptic equations at the boundary of $C^1$ Dini domains, J. Differ. Equations 261(12) (2016), 6718--6757.
\bibitem[DCSV]{DCSV} M. Di Cristo, E. Sincich, S. Vessella, Size estimates of unknown boundaries with Robin type condition, Proceeding of the Royal Society of Edinburgh: Section A, 47 (2017), 727-741.

\bibitem[KN]{l:KuNy} I. Kukavica, K. Nystr\"om, Unique continuation on the boundary for Dini domains, Proc. Amer. Math. Soc. 126(2) (1998), 441--446.

\bibitem[MRVa]{MRV-IUMJ-2007}
A. Morassi, E. Rosset, S. Vessella, Size estimates for inclusions in an elastic plate by boundary measurements, Indiana University Mathematical Journal 56 (2007), 2325--2384.

\bibitem[MRVb]{l:mrv19} A. Morassi, E. Rosset, S. Vessella, Optimal stability in the identification of a
rigid inclusion in an isotropic Kirchhoff-Love plate, SIAM J.
Math. Anal. 51(2) (2019), 731--747.
\bibitem[MRVc]{MRV-LeMatem-2020}
A. Morassi, E. Rosset, S. Vessella, Doubling inequality at the boundary for the Kirchhoff-Love plate's equation with Dirichlet conditions, Le Matematiche LXXV (2020), 27--55.

\bibitem[MRVd]{MRV-JEMS}
A. Morassi, E. Rosset, S. Vessella, Optimal identification of a cavity in the Generalized Plane Stress problem in linear elasticity, J. Eur. Math. Soc. (to appear).

\bibitem [S]{S} 
E. Sincich, Stability for the determination of unknown boundary and impedance with a Robin boundary condition, SIAM J. Math. Anal. 42 (6) (2010), 2922--2943.

\bibitem[TWK]{Timoshenko_Woinowsky-Krieger}
S. Timoshenko, S. Woinowsky-Krieger, \emph{Theory of Plates and Shells}, McGraw-Hill, New York, 1959.



\end{thebibliography}
\end{document}